\documentclass{article}

\usepackage{arxiv}

\usepackage[utf8]{inputenc}
\usepackage[T1]{fontenc}
\usepackage{hyperref}
\usepackage{url}
\usepackage{booktabs}

\usepackage{amsmath,amssymb,eucal}

\usepackage{amsfonts}
\usepackage{graphicx}

\title{Tight bounds on Poisson tails with applications to Szász-Mirakyan operators}

\newif\ifuniqueAffiliation

\uniqueAffiliationtrue

\ifuniqueAffiliation
\author{Jos\'e A. Adell  \\
	Departamento de M\'etodos Estad\'\i sticos\\
	Universidad de Zaragoza\\
	50009 Zaragoza, Spain \\
	\texttt{adell@unizar.es} \\
	\And
	Daniel C\'ardenas-Morales \\
	Departamento de Matem\'aticas\\
	Universidad de Ja\'en\\
	23071 Ja\'en, Spain \\
	\texttt{cardenas@ujaen.es} \\
}
\else

\fi

\renewcommand{\shorttitle}{}

\renewcommand{\headeright}{}

\renewcommand{\undertitle}{}

\hypersetup{pdftitle={},pdfsubject={},pdfauthor={},pdfkeywords={},}

\newtheorem{theorem}{Theorem}
\newtheorem{lemma}[theorem]{Lemma}

\begin{document}

\maketitle

\begin{abstract}
We obtain tight bounds on Poisson tails which are easy to handle. A short proof based on the median of the gamma distribution is given. Numerical comparisons with other known estimates are made. As an application, we consider the rates of pointwise convergence for $S_nf(x)-f(x)$, as $n\rightarrow \infty$, where $S_n$ is the Szász-Mirakyan operator and $x$ belongs to an open interval $I$. If the function $f$ is affine on $I$, the rate of convergence is exponential. This property is no longer true at the boundary of $I$.
\end{abstract}


\section{Introduction and main results}
\label{sec1}

Let $\mathbb{N}$ be the set of positive integers and $\mathbb{N}_0=\mathbb{N}\cup \{0\}$. Throughout this note, we assume that $n\in \mathbb{N}$ and $x>0$. Let $N_{\lambda}$, $\lambda \geq 0$, be a random variable having the Poisson distribution with mean $\lambda$, i.e.,
\begin{equation*}
    P\left(N_{\lambda}=k\right)=\frac{\lambda ^k}{k!}e^{-\lambda},\quad k\in \mathbb{N}_0.
\end{equation*}

We present tight bounds on the probabilities of right and left tail events of the form
\begin{equation}\label{1star}
    P\left(N_{nx}\geq nb\right),\quad  P\left(N_{nx}\leq na\right),\quad 1/n\leq a\leq x \leq b.
\end{equation}
The estimates are given in terms of the non negative convex function
\begin{equation*}
H(t,x)=t \log \frac{t}{x}-t+x, \quad t>0 \qquad (H(0,x)=x),
\end{equation*}
called  the Kullback-Leibler divergence between Poisson random variables with means $t$ and $x$.

Many different upper and lower bounds for the probabilities in (\ref{1star}) are available in the literature. We reparametrize these bounds in order to facilitate their comparison. Boucheron et al. \cite{BoucheronLugosiMassart} (see also \cite{ZhangZhou}) showed that
\[
P\left(N_{nx}\geq nb\right)\leq e^{-nH(b,x)},\quad 0<x\leq b,
\]
as well as
\[
P\left(N_{nx}\leq na\right)\leq e^{-nH(a,x)},\quad 0\leq a\leq x,\quad nx\geq 1.
\]
Denote by $\lfloor \cdot \rfloor$ and $\lceil \cdot \rceil $ the floor and ceiling functions, respectively, and by $\Phi$ the standard normal distribution. Short \cite{Short} showed that for any $c>0$ such that $\lfloor nc \rfloor \geq 1$, we have
\[
\Phi \left( sign (\lfloor nc \rfloor-nx)\sqrt{2nH(\gamma,x)}\right)<P\left(N_{nx}\leq \lfloor nc \rfloor \right)
< \Phi \left( sign (\lfloor nc \rfloor-nx)\sqrt{2nH(\gamma +1/n,x)}\right),\quad  \gamma=\frac{\lfloor nc \rfloor}{n}.
\]
These estimates are derived from an analogous result concerning the binomial distribution proved by Zubkov and Serov \cite{ZubkovSerov}. A more explicit estimate obtained by Short in the same paper is the following. For any $b>0$ such that $\lceil nb \rceil \geq nx+1$, we have
\begin{equation}\label{1starstar}
P\left(N_{nx}\geq \lceil nb \rceil\right)\leq \frac{e^{-nH(\beta -1/n,x)}}{\max \left(2,\sqrt{4 \pi n H(\beta-1/n,x)}\right)},\qquad \beta =\frac{\lceil nb \rceil }{n}.
\end{equation}
A different kind of estimates was proposed by Klar \cite{Klar}. This author showed that for any $b>0$ such that $\lceil nb \rceil +1 > nx$, we have for any $k\in \mathbb{N}$
\[
P\left( \lceil nb \rceil  \leq N_{nx} \leq \lceil nb \rceil +k-1\right) \leq P\left( N_{nx}\geq \lceil nb \rceil \right)
\]
\[
\leq \left(1-\frac{(nx)^k}{\prod_{i=1}^k(\lceil nb \rceil +1)}\right)^{-1}P\left(\lceil nb \rceil \leq N_{nx}\leq \lceil nb \rceil +k-1\right).
\]
More accurate, but also more involved, upper and lower bounds have recently been obtained by From and Swift \cite{FromSwift}.

We state the main results of this note.
\begin{theorem} (Right-tail events)\label{teorema1}
Let $0<x\leq b$ and $\beta=\lceil nb \rceil /n$. Then,
\begin{equation}\label{2}
\frac{\exp (-1/(2 n \beta ))}{\sqrt{2 \pi n \beta }}e^{-nH(\beta,x)}\leq P\left(N_{nx}\geq nb\right)\leq R_n(x,\beta)e^{-nH(\beta ,x)},
\end{equation}
where
\[
R_n(x,\beta)=\min \left(\frac{x}{(\beta -x+1/n)\sqrt{2\pi n \beta}},\frac{1}{2}\right)+\frac{1}{\sqrt{2 \pi n \beta}}.
\]
\end{theorem}

If $b$ is far away from $x$,  then Theorem \ref{teorema1} gives us
\[
P\left(N_{nx}\geq nb\right)\sim \frac{e^{-nH(\beta ,x)}}{\sqrt{2\pi n \beta}}.
\]
On the contrary, if $b\sim x$ then the central limit theorem (see also Lemma \ref{lema3} below) implies that $P\left(N_{nx}\geq nb\right)\sim 1/2$, as $n\rightarrow \infty$. In such a case, the upper bound in (\ref{2}) is approximately $1/2$, since $H(x,x)=0$. We therefore conclude that the upper bound in (\ref{2}) is asymptotically sharp.

\begin{theorem} (Left-tail events)\label{teorema2}
Let $1/n\leq a \leq x$ and $\alpha=\lfloor na \rfloor /n$. Then,
\begin{equation}\label{3}
\frac{\exp (-1/(2 n \alpha  ))}{\sqrt{2 \pi n \alpha }}e^{-nH(\alpha,x)}\leq P\left(N_{nx}\leq na\right)\leq L_n(x,\alpha)e^{-nH(\alpha ,x)},
\end{equation}
where
\[
L_n(x,\alpha)=\min \left(\frac{x}{(x-\alpha +1/n)\sqrt{2\pi n \alpha}},\frac{1}{2}\right)+\frac{1}{\sqrt{2 \pi n \alpha}}.
\]
\end{theorem}

Similar comments to those following Theorem \ref{teorema1} show that the upper bound in (\ref{3}) is asymptotically sharp.

A simple proof of Theorem \ref{teorema1} based on the median of the gamma distribution is given in the following section. A numerical comparative discussion with the bound given in (\ref{1starstar}) is carried out in Section \ref{sec3}. Finally, in Section \ref{sec4}, we apply Theorems \ref{teorema1} and \ref{teorema2} to the Szász-Mirakyan operator $S_n$. In general, the rate of pointwise convergence for $S_nf(x)-f(x)$, as $n\rightarrow \infty$, is $n^{-1}$ at most. However, if $x$ belongs to an open interval $I\subseteq [0,\infty )$ and $f$ is an affine function on $I$, then the rate of convergence becomes exponential. This property is no longer true at the boundary of $I$.  Similar results have recently been achieved for the classical Bernstein polynomials in \cite{AdellCardenasLopez}.

\section{The proofs}
\label{sec2}
Recall that Stirling's approximation (see \cite{Robbins}) states that
\begin{equation}\label{4}
\frac{\exp (-1/(12 k  ))}{\sqrt{2 \pi k }}\leq P\left(N_{k}=k\right)\leq \frac{\exp (-1/(12 k+1 ))}{\sqrt{2 \pi k }},\quad k\in \mathbb{N}.
\end{equation}
On the other hand, the well-known Poisson-Gamma relation says that (cf. Johnson et al. \cite[p. 162]{JohnsonKempKotz})
\begin{equation}\label{5}
P\left(N_{\lambda}\leq k\right)=\int_{\lambda}^{\infty}g_k(u)du,\quad g_k(u)=\frac{u^k}{k!}e^{-u},\quad \lambda > 0,\ k\in \mathbb{N}.
\end{equation}
Finally, the median $\lambda _k$ of the gamma distribution $\Gamma(k+1,1)$ is given by
\begin{equation}\label{6}
\int_0^{\lambda _k}g_k(u)du=\int_{\lambda _k}^{\infty}g_k(u)du=\frac{1}{2},\quad \ k\in \mathbb{N}.
\end{equation}
It was shown in \cite{AdellJodra2005} that
\begin{equation}\label{7}
\lambda _k\in (k,k+1),\quad \ k\in \mathbb{N}.
\end{equation}
Sharp estimates of $\lambda_k$ can be found in \cite{AdellJodra2008}. The following result has interest by itself.

\begin{lemma} \label{lema3}
Let $k\in \mathbb{N}$. Then,
\begin{equation}\label{7star}
\frac{1}{2}  \leq \min \left(P\left(  N_k\geq k\right) ,P\left( N_k\leq k \right)\right)   \leq \max \left(P\left(  N_k\geq k\right) ,P\left( N_k\leq k \right)\right)
 \leq \frac{1}{2}+\frac{1}{\sqrt{2 \pi k}}.
\end{equation}
\end{lemma}
\textit{Proof}. The first inequality in (\ref{7star}) was proved in \cite[Lemma 1]{AdellJodra2005}. Using (\ref{5})-(\ref{7}) and the second inequality in (\ref{4}), we have
\[
P\left( N_k\geq k\right) =P\left( N_k= k\right)+P\left( N_k > k\right)\\ =P\left( N_k =k\right)+\int_0^kg_k(u)du  \leq \frac{1}{\sqrt{2 \pi k}}+\frac{1}{2}.
\]
If $k=1$, then
\[
P\left( N_1 \leq 1\right)=\frac{2}{e}\leq \frac{1}{2}+\frac{1}{\sqrt{2 \pi}}.
\]
If $k\geq 2$, we have as above
\[
P\left( N_k\leq k\right) =P\left( N_k= k\right)+P\left( N_k \leq  k-1\right)=P\left( N_k =k\right)+\int_k^{\infty }g_{k-1}(u)du  \leq \frac{1}{\sqrt{2 \pi k}}+\frac{1}{2}.
\]
This completes the proof.
\hfill $\Box$

\subsection{Proof of Theorem \ref{teorema1}}
\label{subsec1}
Let $k\in \mathbb{N}_0$. We start with the basic identity
\begin{equation}\label{8}
P\left( N_{nx}=n\beta +k\right)=P\left( N_{n\beta }=n \beta +k\right)\left(\frac{x}{\beta}\right)^ke^{-nH(\beta, x)}.
\end{equation}
Observe that
\[
\frac{P\left( N_{n\beta }=n\beta +k\right)}{P\left( N_{n\beta }=n\beta \right)}\leq \left(\frac{n\beta }{n\beta +1}\right)^k,\quad k\in \mathbb{N}.
\]
Hence, we have from (\ref{8}) and the second inequality in (\ref{4})
\begin{align}\label{9}
P\left( N_{nx} \geq n\beta \right)  \leq P\left( N_{n\beta }=n\beta \right)e^{-nH(\beta ,x)}\left( 1+\sum_{k=1}^{\infty}\left( \frac{nx}{n\beta +1}\right)^k\right) 
 \leq \left( 1+\frac{x}{\beta -x+1/n} \right)\frac{e^{-nH(\beta ,x)}}{\sqrt{2 \pi n \beta }}.
\end{align}
Since $x\leq \beta$, we get from (\ref{8}) and Lemma \ref{lema3}
\[
P\left( N_{nx} \geq n\beta \right) \leq P\left( N_{n\beta } \geq n\beta \right)e^{-nH(\beta ,x)} \leq \left( \frac{1}{2}+\frac{1}{\sqrt{2 \pi n \beta }}\right)e^{-nH(\beta,x)}.
\]
This and (\ref{9}) show the upper bound in (\ref{2}). With respect to the lower bound, we use again (\ref{8}) and the first inequality in (\ref{4}) to obtain
\[
P\left( N_{nx} \geq n\beta \right) \geq P\left( N_{n\beta } = n\beta \right)e^{-nH(\beta ,x)} \geq \left( \frac{\exp \left( -1/(12 n \beta )\right)}{\sqrt{2 \pi n \beta }}\right)e^{-nH(\beta,x)}.
\]
The proof is complete.

\subsection{Proof of Theorem \ref{teorema2}}
\label{subsec2}
The proof follows along the lines of that of Theorem \ref{teorema1} by using the identity
\[
P\left( N_{nx}=n\alpha -k\right) = P\left( N_{n\alpha}=n\alpha -k\right) \left( \frac{\alpha}{x}\right)^ke^{-nH(\alpha ,x)},\quad k=0,1,\ldots, n\alpha .
\]
and therefore we omit it.

\section{Numerical comparative discussion}
\label{sec3}
For the sake of brevity, we only compare the upper bound in Theorem \ref{teorema1} with estimate (\ref{1starstar}) obtained by Short \cite{Short}. Let $b>0$ and $n\in \mathbb{N}$ be fixed. A fortiory, $\beta =\lceil nb\rceil /n $ is also fixed.

Denote
\[
x_{\beta,n}=\frac{(\beta +1/n)\sqrt{2 \pi n \beta}}{2+\sqrt{2 \pi n \beta}}.
\]
Observe that the upper bound in Theorem \ref{teorema1} can be written as
\begin{equation}\label{A}
P\left( N_{nx} \geq nb\right)\leq
\left\{
  \begin{array}{ll}
    \left( \frac{x}{\beta -x+1/n}+1 \right)\frac{\exp (-nH(\beta,x))}{\sqrt{2 \pi n \beta}}, & 0<x\leq x_{\beta,n} \\
     &  \\
   \left( \frac{1}{2}+\frac{1}{\sqrt{2 \pi n \beta}}\right)e^{-nH(\beta ,x)}, & x_{\beta,n} <x \leq b.
  \end{array}
\right.
\end{equation}

On the other hand, the function
\[
g_{\beta,n}(x)=\pi nH(\beta -1/n,x),\quad 0<x\leq \beta -1/n,
\]
strictly decreases from $\infty $ to $0$. Thus, there exists a unique $y_{\beta,n}\in (0,\beta -1/n)$ such that $g_{\beta,n}(y_{\beta,n})=1$. Observe that (\ref{1starstar}) can be written as
\begin{equation}\label{1starstarstar}
P\left(N_{nx}\geq \lceil nb \rceil\right)\leq
\left\{
  \begin{array}{ll}
  \frac{\exp (-nH(\beta - 1/n,x))}{\sqrt{4 \pi n H(\beta -1/n,x)}}, & 0<x\leq y_{\beta,n} \\
     &  \\
    \frac{\exp (-nH(\beta -1/n,x))}{2}, & y_{\beta,n} <x \leq \beta -\frac{1}{n}.
  \end{array}
\right.
\end{equation}

Note that (\ref{A}) holds for $x\leq b$, while (\ref{1starstarstar}) holds for $x\leq \beta-1/n$. Hence, to compare both estimates, we consider only the case $x\in (0,\beta -1/n)$. Besides, we assume that $n\beta =\lceil nb\rceil \geq 9$. It can be proved by analytic methods that under this mild assumption inequalities $y_{\beta,n}<x_{\beta,n}<\beta -1/n$ are fulfilled.

\begin{figure}[h!]
    \centering
\includegraphics[width=0.6\textwidth]{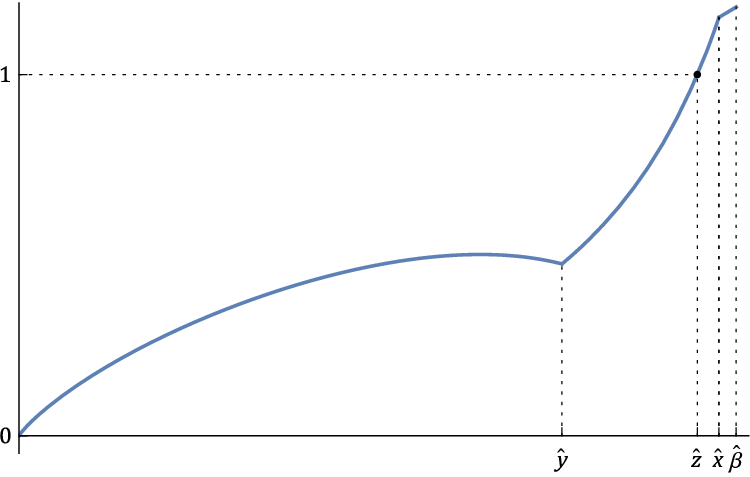}
    \caption{Graph of $Q_{\beta,n}(x)$. $\hat{x}=x_{\beta,n}$. $\hat{y}=y_{\beta,n}$. $\hat{z}=z_{\beta,n}$. $\hat{\beta}=\beta -1/n$.}
    \label{figura1}
\end{figure}

Figure \ref{figura1}, plotted only for illustrative purposes, shows the graph of the quotient $Q_{\beta,n}(x)$, $0<x\leq \beta -1/n$, between the upper bounds in (\ref{A}) and (\ref{1starstarstar}). An elementary, but also cumbersome, analysis shows  that $Q_{\beta,n}(x)$ is continuous, less than $1$ on $(0,y_{\beta,n})$, strictly increasing on $(y_{\beta,n},x_{\beta,n})$ and $(x_{\beta,n},\beta-1/n)$, and takes the value $1$ at a point $z_{\beta,n}\in (y_{\beta,n},x_{\beta,n})$, whose explicit expression in terms of $\beta$ and $n$ turns out to be
\[
z_{\beta,n}=\frac{\beta +1/n}{1+e\sqrt{2/\pi}(n\beta +1)(n\beta )^{-n\beta -1/2}(n\beta -1)^{n\beta -1}}.
\]
Therefore, the estimate in (\ref{A}) is better than that in (\ref{1starstarstar}) in the interval $(0,z_{\beta,n})$.

\section{The Szász-Mirakyan operator}
\label{sec4}
The Szász-Mirakyan operator is defined by
\[
S_nf(x)=\sum_{k=0}^{\infty}f\left(\frac{k}{n}\right)\frac{(nx)^k}{k!}e^{-nx}=\mathbb{E}f\left( \frac{N_{nx}}{n}\right),\quad x\geq 0,
\]
where $f:[0,\infty )\rightarrow \mathbb{R} $ is any arbitrary function for which the preceding series are convergent for any $x\geq 0$. With regard to the approximation properties of this operator, consider the Ditzian-Totik modulus of smoothness of $f$ with step weight function $\varphi (x)=\sqrt{x}$, $x\geq 0$, that is,
\[
\omega_2^{\varphi}(f;\delta)=\sup \{|f(x-h\sqrt{x})-2f(x)+f(x+h\sqrt{x})|: x\geq h^2, 0\leq h\leq \delta            \},\ \delta \geq 0.
\]
Denote by $\| \cdot \|$ the usual supremum norm on $[0,\infty )$. It is known (cf. Totik \cite{Totik}) that if $f:[0,\infty )\rightarrow \infty $ is a continuous function such that $\omega_2^{\varphi} (f;1)<\infty $, then
\begin{equation}\label{10}
K_1\omega_2^{\varphi}\left(f;\frac{1}{\sqrt{n}}\right)\leq \|S_nf-f\| \leq K_2\omega_2^{\varphi}\left(f;\frac{1}{\sqrt{n}}\right ),
\end{equation}
for some absolute positive constants $K_1$ and $K_2$.
P\u{a}lt\u{a}nea \cite{Paltanea} proved that (\ref{10}) is fulfilled for $K_2=5$, whereas in \cite{AdellCardenas} the authors showed that we can take $K_2=2.43$. However, it seems that no specific value for the lower constant $K_1$ has been provided yet.

Note that the rate of convergence in (\ref{10}) is $n^{-1}$ at most. This rate of convergence can be improved if the function $f$ is smooth. In fact, López-Moreno and Latorre-Palacios \cite{LopezLatorre} showed that if $f^{(j)}(x_0)=0$, $j\in \mathbb{N}_0$, for some $x_0>0$, then $S_nf(x_0)=o(n^{-k})$, for all $k\in \mathbb{N}$.

Denote by $1_A$ the indicator function of the set $A$. For locally affine functions $f$ satisfying appropriate integrability conditions, the rate of pointwise convergence becomes exponential, as shown in the following result.

\begin{theorem}\label{teorema3}
Let $x\in (a,b)$, with $1/n\leq a<b<\infty $. Suppose that $f,\ell :[0,\infty )\rightarrow \mathbb{R}$ are functions such that $f$ is arbitrary, $\ell$ is affine, and $f(t)=\ell (t)$, $t\in (a,b)$. If \begin{equation}\label{12}
\mathbb{E}|f|^p\left(\frac{N_{nx}}{n}\right)<\infty ,
\end{equation}
for some $p>1$, then
\[
|S_nf(x)-f(x)|\leq \mathbb{E}^{1/p}|f-\ell |^p\left( \frac{N_{nx}}{n}\right)\left( \frac{L_n^{1/q}(x,\alpha)}{e^{nH(\alpha,x)/q}}+ \frac{R_n^{1/q}(x,\beta)}{e^{nH(\beta,x)/q}}\right),
\]
where $p^{-1}+q^{-1}=1$.
\end{theorem}
\textit{Proof}. Consider the function $g(t)=f(t)-\ell (t)$, $t\geq 0$. Since $g(x)=0$ and $S_n\ell (x)=\ell (x)$, we have
\begin{align*}
|S_nf(x)-f(x)| =|S_n(g+\ell )(x)-(g+\ell )(x)|=|S_ng(x)| \\
\leq \mathbb{E}|g|\left(   \frac{N_{nx}}{n}  \right) 1_{\{ N_{nx}\leq na \}}
+\mathbb{E}|g|\left( \frac{N_{nx}}{n} \right) 1_{\{ N_{nx}\geq nb \}}.
\end{align*}
Hence, the result follows from (\ref{2}), (\ref{3}), and Hölder's inequality.
\hfill $\Box$

\noindent \textbf{Remark 5} 
If, instead of condition (\ref{12}), we assume that $\|f-\ell \|<\infty$, then the same proof as that of Theorem \ref{teorema3} gives us
\[
|S_nf(x)-f(x)|\leq \|f-\ell \|\left( L_n(x,\alpha )e^{-nH(\alpha,x)}+R_n(x,\beta)e^{-nH(\beta,x)} \right).
\]
If $b=\infty$, we similarly have
\[
|S_nf(x)-f(x)|\leq \|f-\ell \|L_n(x,\alpha )e^{-nH(\alpha,x)},\quad 1/n\leq a<x.
\]

\ 

We finally point out that Theorem \ref{teorema3} and Remark 5 may fail at the boundary of $(a,b)$. Indeed, consider the function
\[
f_s(t)=(t-a)_-^s,\quad t\geq 0,\ s>0,\ 1/n\leq a,
\]
where $y_-=\max (0,-y)$, $y\in \mathbb{R}$. Using the convergence of moments in the central limit theorem (see, for instance, Billingsley \cite[p. 338]{Billingsley}), we get
\[
S_nf_s(a)-f_s(a)=\left(\frac{a}{n}\right)^{s/2}\mathbb{E}\left(\frac{N_{na}-na}{\sqrt{na}}\right)_-^s\sim \left(\frac{a}{n}\right)^{s/2}\mathbb{E}Z_-^s,\quad  n\rightarrow \infty ,
\]
where $Z$ is a random variable having the standard normal distribution.

\

\textbf{Corresponding author}

Correspondence to Daniel Cárdenas-Morales.

\noindent \textbf{Funding}

The first author is supported by Research Project DGA (E48\_23R). The second author is supported by Junta de Andaluc\'\i a (Research Group FQM-0178).

\noindent \textbf{Declaration of competing interest}

The authors declare that they have no competing interests.

\noindent\textbf{Data availability}

Data were not generated or analyzed for the research described in the article.

\bibliographystyle{unsrtnat}

\end{document}